\documentclass{amsart}

\usepackage[T1]{fontenc}
\usepackage{textcomp}
\usepackage{amsfonts}
\usepackage{amsthm}
\usepackage{ amssymb }
\usepackage{pst-node}
\usepackage{tikz-cd}
\usetikzlibrary {arrows.meta,automata,positioning}
\usetikzlibrary{positioning,decorations.pathreplacing}
\usepackage{ bbold }
\usepackage[all]{xy}
\usepackage{mathtools}
\usepackage{chngcntr}
\usepackage{thmtools}
\usepackage[acronym,nogroupskip,nonumberlist,nopostdot,toc]{glossaries}
\usepackage{glossaries}
\usepackage{mathrsfs}
\usepackage{caption}
\usepackage{subcaption}
\usepackage{tcolorbox}
\usepackage{listings}
\usepackage{imakeidx}
\usepackage{enumerate}
\usepackage{enumitem}
\usepackage{hyperref}
\usepackage[nameinlink, capitalize, noabbrev]{cleveref}

\hypersetup{ colorlinks=true,linkcolor=teal,    citecolor=magenta, }

\def\set#1{{\def\st{\;:\;}\left\{#1\right\}}}
\def\abs#1{\left\vert{#1}\right\vert}
\def \<#1>{{\left\langle{#1}\right\rangle}}

\def\ZZ{\mathbb Z}

\def\OK{\mathcal{O}_K}

\def\NN{\mathbb N}

\def\Q-{\overline{\mathbb Q}}

\DeclareMathOperator{\Aut}{Aut}

\DeclareMathOperator{\IMG}{IMG}

\DeclareMathOperator{\St}{St}
\DeclareMathOperator{\RiSt}{RiSt}
\DeclareMathOperator{\rist}{rist}
\DeclareMathOperator{\Sym}{Sym}

\DeclareMathOperator{\st}{st}

\DeclareMathOperator{\ab}{ab}

\DeclareMathOperator{\SL}{SL}

\newtheorem{Theorem}{Theorem}[section]

\newtheorem{Proposition}[Theorem]{Proposition}
\newtheorem{Lemma}[Theorem]{Lemma}

\newtheorem{theorem}{Theorem}[section]

\title[A just-infinite iterated monodromy group without CSP]{A just-infinite iterated monodromy group without the congruence subgroup property}
\author{Santiago Radi}
\address{Santiago Radi: Department of Mathematics, Texas A\&M University, 77843 College Station, U.S.A.
}
\email{santiradi@tamu.edu}
\thanks{The author is supported by Grigorchuk's Simons Foundation Grant MP-TSM-00002045 and the department of Mathematics of Texas A\&M University.}
\date{May 2025}

\begin{document}

\begin{abstract}
We prove that the iterated monodromy group of the polynomial $z^2+i$ is just-infinite, regular branch and does not have the congruence subgroup property. This yields the first example of an iterated monodromy group of a polynomial with these properties. Additional information is provided about the congruence kernel, rigid kernel and branch kernel of this group.
\end{abstract}

\maketitle

\section{Introduction}
Given a residually finite group $G$, one can fix a sequence of epimorphisms $\set{\varphi_n}_{n \geq 0}$ from $G$ onto finite groups such that the sequence of normal subgroups  $\set{\ker(\varphi_n)}_{n \geq 0}$ is partially ordered with respect to the inclusion and we also have $\bigcap_{n \geq 0} \ker(\varphi_n) = 1$. Each subgroup in this sequence is called a principal congruence subgroup, and a subgroup of $G$ is called a congruence subgroup if it contains a principal congruence subgroup. 

The congruence subgroup problem asks whether every finite-index subgroup in $G$ is a congruence subgroup. If this happens, we say that $G$ has the congruence subgroup property with respect to a given sequence. For example, it is proved by Bass, Milnor and Serre that if $K$ is a number field and $\OK$ is its ring of integers, then $\SL_n(\OK)$ has the congruence subgroup property whenever $n \geq 3$ and it does not have the congruence subgroup property if $n = 2$ (see \cite{BassMilnorSerre1967}).

As observed in \cite{GrigorchukNekrashevichSuschanski2000}, every countable residually finite group can be embedded in the group of automorphisms of a spherically homogeneous infinite rooted tree $T$, which is denoted $\Aut(T)$. The group $\Aut(T)$ is residually finite and the principal congruence subgroups used are the stabilizers of level $n$ for each $n \geq 0$. The stabilizer of level $n$, denoted $\St(n)$, consists of all the elements that fix all the vertices at the level $n$ of the tree. If $G$ is a subgroup of $\Aut(T)$, we can consider $\set{\St(n) \cap G}_{n \geq 0}$ as the principal congruence subgroups of $G$. 

The study of the congruence subgroup property in groups acting on rooted trees (more specifically for branch groups) was initiated by Grigorchuk, Herfort and Zalesskii in \cite{GrigHerfortZalesskii1998} and since then, few examples of groups with or without the congruence subgroup property are known in the literature. In \cite{Grigorchuk2000}, Grigorchuk proved that the group nowadays known as the first Grigorchuk group (see \cite{Grig1980}) has the congruence subgroup property (see \cite[Proposition 10]{Grigorchuk2000}). In \cite{FernandezAlcoberGarrido2017}, Fernández-Alcober, Garrido and Uria-Albizuri proved that the class of Grigorchuk-Gupta-Sidki (GGS) groups (see \cite{GGS1983}) with non-constant defining vector also have the congruence subgroup property. On the other hand, Pervova constructed the first family of groups without the congruence subgroup property in \cite{Pervova2007}, by defining the groups EGS, a slightly different version of the GGS groups. Generalizing the result of Pervova, Thillaisundaram and Uria-Albizuri in \cite{ThillaisundaramUria2021} studied the multi-EGS groups (defined in \cite{KlopschThillaisundaram2018}), classifying the cases where they do not have the congruence subgroup property. In \cite{Bartholdi2012}, Bartholdi, Siegenthaler and Zalesskii proved that the twisted twin of the Grigorchuk group and the Hanoi tower group (defined in \cite{GrigNekraSunic2006}) do not have the congruence subgroup property. Finally, in \cite{Garrido2016}, Garrido proved that in the case of branch groups, the congruence subgroup property is independent of the branch action of the group.

Within the groups acting on trees, one important family is the family of iterated monodromy groups, used by Bartholdi and Nekrashevych to solve the Twisted Rabbit Problem in \cite{NekraBartholdi2005}. One celebrated example, and the main object of this article, is the iterated monodromy group of the polynomial $z^2+i$. This 3-generated group that acts on the binary tree was initially considered in \cite{BartholdiGrigNekra2003} and it was proved in \cite{BuxPerez2006} that the group has intermediate growth. In \cite{GrigSavchukSunic2007}, authors proved that this group is regular branch (a particular case of branch group) over a maximal branching subgroup which we will denote $K$ in this article, and found an $L$-presentation, i.e. a presentation of a group by generators and relations that involves a finite set of relators and their iterations by a substitution (see \cite{Self_similar_groups}). In \cite{GrigSavchukSunic2007}, it was also shown that the group exhibits interesting spectral properties. Recently, the authors in \cite{Savchukcryptography2024} study applications of this group (and other contracting groups) as a group-based cryptosystem. 

Our main result is the following:

\begin{theorem}
The group $\IMG(z^2+i)$ is just-infinite and does not have the congruence subgroup property.
\label{theorem: IMG joo and no CSP}
\end{theorem}

Recall that a group is just-infinite if every non-trivial quotient is finite. This constitutes the first example answering the congruence subgroup property of a just-infinite iterated monodromy group of a quadratic complex polynomial.

The properties just-infinite and congruence subgroup property are related in the case of finitely generated groups acting on trees. It was proved by Francoeur in \cite[Theorem 4.4.4]{Francouer2019} that groups with the congruence subgroup property are just-infinite. In virtue of our result in \cref{theorem: IMG joo and no CSP}, the group $\IMG(z^2+i)$ constitutes a counter example to the converse of \cite[Theorem 4.4.4]{Francouer2019}. This group is not the first counterexample of this kind, as the majority of the branch groups without the congruence subgroup property previously mentioned are also just-infinite.

The strategy in this article will be based on the maximal regular branching subgroup $K$ and the $L$-presentation of $\IMG(z^2+i)$ found in \cite{GrigSavchukSunic2007}. If we let $C_m$ denote the cyclic group with order $m$, the steps are the following:

\begin{enumerate}
\item We find $S_K = \set{x,y,z,w,t}$, a set of five generators of $K$ where all the elements have order $4$. Therefore, we have $K/K' \leq C_4^5$. As $\IMG(z^2+i)$ is regular branch over $K$, this in particular implies that $\IMG(z^2+i)$ is just-infinite by \cite[Theorem 4]{Grigorchuk2000}. In order to prove that $\IMG(z^2+i)$ does not have the congruence subgroup property, the goal will be to prove that $K/K' \simeq C_4^5$, namely, the quotient has no other relations.

\item We consider $\bigoplus_{\NN} S_K$, the free monoid over the alphabet $S_K$ and we define the morphism of monoids $e: \bigoplus_{\NN} S_K \rightarrow \ZZ/2\ZZ$, that returns the parity of the length of the word. Using the $L$-presentation of $\IMG(z^2+i)$ and the Reidemeister algorithm, we find an $L$-presentation for $K$. With this presentation, we prove that the function $e$ descends to a function on $K$.

\item Using the action of $K$ on the first $n$ levels of the tree, we deduce that the only possible relations are $ty^{-2} \in K'$ and $wz^{-2} \in K'$.

\item Using the function $e$ and a subgroup $\mathcal{K}$, which depends on $e$ and satisfies $K' \leq \mathcal{K} \leq K$, we deduce that $ty^{-2}$ and $wz^{-2}$ are not in $\mathcal{K}$.
\end{enumerate}

In conclusion:

\begin{theorem}
Let $K$ be the maximal branching subgroup of $\IMG(z^2+i)$. Then $K/K'$ is isomorphic to $C_4^5$.
\label{theorem: K/K' is C45}
\end{theorem}

Although it is not necessary for the proof of \cref{theorem: K/K' is C45}, we will prove that $K/\mathcal{K} \simeq C_4 \times C_2 \times C_2$. As a consequence, the quotient $\mathcal{K}/K'$ is an abelian group of order $64$.

Using the algorithm presented in \cite[Remark 2.8]{Bartholdi2012}, we also obtain:

\begin{theorem}
The rigid kernel of $\IMG(z^2+i)$ is trivial whereas the congruence kernel and the branch kernel are isomorphic to $C_4[[\partial T]]$, namely, the abelian group of formal sums of elements in $\partial T$ with coefficients in $C_4$, where $\partial T$ represents the boundary of the binary tree $T$.  
\label{theorem: profinite kernels IMG(z2+i)}
\end{theorem}

For definitions of rigid kernel and congruence kernel, see \cref{section: congruence subgroup property}.

\subsection*{Organization}
The article is organized in seven sections. In \cref{section: preliminaries}, the background necessary to follow this article is introduced. \cref{section: The group IMG(z2+i)} reviews known results about the group $\IMG(z^2+i)$ that will be used in the subsequent sections. \cref{section: generators and presentation of the branching subgroup} is focused on point (1) of the strategy, where we determine two sets of generators and an $L$-presentation for $K$. \cref{section: The exponent function in K} corresponds to point (2) of the strategy, proving that the function $e$ descends to a well-defined function on $K$. In \cref{section: abelianization of K}, \cref{theorem: IMG joo and no CSP} and \cref{theorem: K/K' is C45} are proved. Finally, in \cref{section: branch kernel and rigid kernel}, the rigid, branch and congruence kernels are calculated, proving \cref{theorem: profinite kernels IMG(z2+i)}.

\subsection*{Acknowledgments}

The author would like to thank R. Grigorchuk, V. Nekrashevych and J. Fariña-Asategui for the valuable discussions held and remarks. Also to the anonymous reviewer for valuable feedback and alternative proofs of some results in \cref{section: abelianization of K}. 

\section{Preliminaries}
\label{section: preliminaries}

\subsection{About general notation} 

We write $H\leq G$ (resp. $H\leq_fG$) to indicate that $H$ is a subgroup of $G$ (resp. $H$ is a subgroup of finite index). The commutator of elements $x,y$ is $[x,y] = xyx^{-1}y^{-1}$. The commutator subgroup of $G$ is denoted by $G'$, and the abelianization of $G$ by $G^{ab} := G/G'$. Given $x,y,z \in G$, we denote $[x,y,z] = [[x,y],z]$. If $K$ and $H$ are subgroups of $G$, then $K^H$ is defined as the group generated by elements of the form $hkh^{-1}$ with $k \in K$ and $h \in H$. We write $C_m$ for a cyclic group of order $m$.

\subsection{Basic lemmas in group theory}

In this subsection we will list basic results that we will use throughout the article.

\begin{Lemma}[{\cite[Proposition 5.1.7]{Robinson1996}}]
\label{lemma: commutator generators}
Let $G$ be a group, $K$ a subgroup of $G$ and $S, S'$ two generating sets of $K$. Then, the commutator subgroup of $K$ is given by $$K' = \<[s,s']:s \in S,s' \in S'>^K.$$ Moreover, if $K \lhd G$ and $H$ is any subgroup containing $K$, then $$K' = \<[s,s']:s \in S,s' \in S'>^H.$$ 
\end{Lemma}

\begin{Lemma}
Let $A$ be a group, $\pi: A \rightarrow \pi(A)$ a surjective map, $B \lhd A$ and call $p_1: A \rightarrow A/B$ and $p_2: \pi(A) \rightarrow \pi(A)/\pi(B)$ the natural projections. Then we have a commutative diagram

\begin{equation*}
\centering
\begin{tikzcd}
A \arrow{rr}{p_1} \arrow{dd}{\pi} & & A/B \arrow{dd}{\varphi} \\
 & & \\
\pi(A) \arrow{rr}{p_2} & & \pi(A)/\pi(B) \\
\end{tikzcd}
\end{equation*}
where $\varphi$ is a well-defined surjective map and $\ker(\varphi) \simeq \ker(\pi)/(\ker(\pi) \cap B)$. Moreover, we have the following equivalences:

\begin{enumerate}
\item $\varphi$ is an isomorphism,
\item $\ker(\pi) \leq B$,
\item $[A:B] = [\pi(A):\pi(B)]$. 
\end{enumerate}
\label{lemma: map defined in abelianizations}
\end{Lemma}

\begin{proof}
First, since $\pi$ is surjective and $B \lhd A$, then $\pi(B) \lhd \pi(A)$, ensuring that $\varphi$ is well defined. Moreover, since all the maps involved are surjective, $\varphi$ is also surjective.

For the kernel, we have the following equivalences:
\begin{align*}
g B \in \ker(\varphi) \Longleftrightarrow \varphi(gB) = p_2 \, \pi(g) = 1 \Longleftrightarrow  \pi(g) \in \pi(B) \Longleftrightarrow \\ g \in \pi^{-1}(\pi(B)) = \<\ker(\pi),B> = \ker(\pi)B.
\end{align*}
Hence $$\ker(\varphi) = \ker(\pi)B/ B \simeq \ker(\pi)/(\ker(\pi) \cap B)$$ by the second isomorphism theorem.

For the equivalences, the first two follow from the form of the kernel of $\varphi$. Indeed, $$\ker(\varphi) = 1 \Longleftrightarrow \ker(\pi) \cap B = \ker(\pi) \Longleftrightarrow \ker(\pi) \leq B.$$

Finally, since we already know $\varphi$ is a surjective morphism of groups, 
\begin{equation*}
\varphi \text{ is an isomorphism } \Longleftrightarrow \abs{A/B} = \abs{\pi(A)/\pi(B)}. \qedhere
\end{equation*}
\end{proof}

\subsection{Groups acting on rooted trees}
\label{section: Groups acting on trees}

A \textit{d-regular rooted tree} $T$ is a tree with infinitely many vertices and a root $\emptyset$, where all the vertices have exactly $d$ descendants. In this article we will in general work with the binary tree, so $d = 2$. The binary tree may also be identified with the finite words over the set $X = \set{1, 2}$. The set of vertices at distance exactly $n\geq 0$ from the root forms the $n$-th level of $T$, denoted $\mathcal{L}_n$. The set of vertices whose distance is at most $n$ from the root forms the $n$-th truncated tree of $T$, denoted by $T^n$. The set of infinite paths starting from the root, or equivalently, the infinite words in $X$, is called the boundary of the tree and is denoted by $\partial T$. 

The group of automorphisms of $T$, denoted $\Aut(T)$, is the group of bijective functions from $T$ to $T$ that preserve adjacency between vertices. This in particular implies that $\Aut(T) \curvearrowright \mathcal{L}_n$ for all $n \in \NN$. For any vertex $v \in T$ the subtree rooted at $v$, which is again a binary tree, is denoted $T_v$. The action of $\Aut(T)$ on $T$ will be on the left, meaning that if $g,h \in \Aut(T)$ and $v \in T$, then $(gh)(v) = g(h(v))$. 

Given a vertex $v \in T$, we denote by $\st(v)$ the \textit{stabilizer of the vertex $v$}, namely, the subgroup of the elements $g \in \Aut(T)$ such that $g(v) = v$. Given $n \in \NN$, we write $\St(n) = \bigcap_{v\in \mathcal{L}_n}\st(v)$ for the \textit{stabilizer of level $n$}. It is not hard to see that $\St(n)$ is a normal subgroup of finite index. We endow $\Aut(T)$ with a topological group structure by declaring $\set{\St(n)}_{n \geq 0}$ to be a base of neighborhoods of the identity. This makes $\Aut(T)$ a profinite group, homeomorphic and isomorphic to the inverse limit $$\varprojlim \Aut(T)/\St(n).$$ The transition maps are given by $\pi_n: \Aut(T) \rightarrow \Aut(T^n)$, where $\pi_n$ restricts the action of an element of $\Aut(T)$ to the first $n$ levels of the tree. 

Let $v \in T$ and $g \in \Aut(T)$. By preservation of adjacency, there exists a unique map $g|_v \in \Aut(T_v)$ such that for all $w \in T_v$, we have $$g(vw)=g(v)g|_v(w).$$ This map is called the \textit{section} of $g$ at the vertex $v$. The sections satisfy the following two fundamental properties: 
\begin{align}
(gh)|_v = g|_{h(v)} h|_v \text{ and } (g^{-1})|_v = (g|_{g^{-1}(v)})^{-1}.
\label{equation: properties sections}
\end{align}
for all $g,h \in\Aut(T)$ and $v \in T$.

The properties in \cref{equation: properties sections} yield the following isomorphism for any natural number $n \in \NN$: 
\begin{align}
\begin{split}
\psi_n:\Aut(T) \to (\Aut(T_{v_1}) \times \dotsb \times \Aut(T_{v_{N_n}}) \rtimes \Aut(T^n) \\ \text{ via } g \mapsto (g|_{v_1},\dotsc,g|_{v_{N_n}}) \pi_n(g)
\end{split}
\label{equation: Aut and semidirect product}
\end{align}
where $v_1,\dotsc,v_{N_n}$ are all the distinct vertices in $\mathcal{L}_n$.

To describe the elements of $\Aut(T)$, we will use the isomorphism from \cref{equation: Aut and semidirect product} with $n = 1$. In some cases, we will also use the notation $(g_0,\dots,g_N)_n$ that corresponds to the image of an element in $\St(n)$ via the map $\psi_n$. 

Now, let us fix a subgroup $G \leq \Aut(T)$. The group $G$ inherits natural actions on $T$ and each set $\mathcal{L}_n$. We define the vertex stabilizers and level stabilizers of $G$ by restricting the corresponding stabilizers of $\Aut(T)$ to $G$, i.e, $$\st_G(v) := \st(v) \cap G \text{ and } \St_G(n) := \St(n) \cap G$$ for $v \in T$ and $n \in \NN$. 

A group $G \leq \Aut(T)$ is \textit{level-transitive} if the actions $G \curvearrowright \mathcal{L}_n$ are transitive for all $n \in \NN$. The group $G$ is said to be \textit{self-similar} if $g|_v \in G$ for all $v \in T$ and $g \in G$. 

Given $G \leq \Aut(T)$, a vertex $v \in T$, and $n \in \NN$, we define the \textit{rigid vertex stabilizer} of $v$, denoted by $\rist_G(v)$, as the subgroup of $G$ consisting of automorphisms that fix every vertex outside the subtree $T_v$. The \textit{$n$-th rigid level stabilizer}, denoted $\RiSt_G(n)$, is the subgroup generated by the subgroups $\rist_G(v)$ with $v \in \mathcal{L}_n$. It is not hard to see that $\RiSt_G(n) \lhd G$ and that $\RiSt_G(n) \leq \St_G(n)$. 

We say that a group $G \leq \Aut(T)$ is a \textit{branch group} if $\RiSt_G(n) \leq_f G$ for all $n \in \NN$.

\subsection{Regular branch groups}

Within the class of branch groups, we have the class of regular branch groups. 

If $K$ is a subgroup of $\Aut(T)$ and $n \in \NN$, we define the \textit{geometric product} of $K$ on level $n$ as $$K_n := \set{g \in \St(n): g|_v \in K \text{ for all $v \in \mathcal{L}_n$}}.$$

Using the isomorphism in \cref{equation: Aut and semidirect product}, if $g \in K_n$, we can represent it as $$g = (k_1, \dots, k_N)_n,$$ where $N = \# \mathcal{L}_n = 2^n$ and $k_i \in K$ for all $i = 1, \dots, N$. Since $K$ is a subgroup, then each element $g_i := (1, \dots, 1, k_i, 1, \dots, 1)_n \in K_n$ and so $g$ can be split as the product of the elements $g_i$. We conclude that if $S$ is a generating set for $K$, then 
\begin{equation}
\set{ (1, \dots, 1, s, 1, \dots, 1)_n: s \in S, i = 1, \dots, 2^n}
\label{equation: generators Kn}
\end{equation}
is a generating set for $K_n$. 

A subgroup $G \leq \Aut(T)$ is called \textit{regular branch} over a subgroup $K$ if $$K_1 \leq_f K \leq_f G.$$ Notice that by iterating the inclusion $K_1 \leq K$, it follows that $K_n \leq K$ for all $n \in \NN$. 

\begin{Lemma}[{\cite[Lemma 10]{Sunic2006}}] 
If $G \leq \Aut(T)$ is a self-similar regular branch group branching over a subgroup $K$ containing $\St_G(m)$ for some $m \in \NN$, then for all $n \geq m$, $$\St_G(n) \simeq \prod_{i = 1}^{2^{n-m}} \St_G(m)$$ via the map $\psi_{n-m}$ defined in \cref{equation: Aut and semidirect product}, restricted to $\St_G(n)$.
\label{lemma: St(n+1) = prod St(n)}
\end{Lemma}

The following lemma will be useful in the proof of \cref{theorem: IMG joo and no CSP}:

\begin{Lemma}
\label{lemma: stabilization pin(K)ab}
Let $T$ be a $d$-regular rooted tree and let $G \leq \Aut(T)$ be a regular branch subgroup branching over $K$. Suppose that $K$ contains $\St_G(m)$ for some $m \in \NN$ and that there exists $n_0 \geq m$ such that $\pi_{n_0+1}(\St_G(n_0)) \leq \pi_{n_0+1}(K')$. Then for all $n \geq n_0$, $$\pi_n(K)/\pi_n(K') \simeq \pi_{n_0}(K)/\pi_{n_0}(K').$$  
\end{Lemma}

\begin{proof}
It is enough to prove that $$\pi_{n+1}(K)/\pi_{n+1}(K') \simeq \pi_n(K)/\pi_n(K')$$ for all $n \geq n_0$. We take $n \geq n_0$ and apply \cref{lemma: map defined in abelianizations} with $A = \pi_{n+1}(K)$, $B = \pi_n(K')$ and $\pi$ the natural restriction from level $n+1$ to level $n$. Then, the induced map $\varphi: \pi_{n+1}(K)/\pi_{n+1}(K') \rightarrow \pi_n(K)/\pi_n(K')$ is an isomorphism if and only if $$\ker(\pi) = \pi_{n+1}(\St_G(n)) \leq \pi_{n+1}(K').$$ Thus, our goal is to prove this inclusion for all $n \geq n_0$.

First of all, we have $(K')_n \leq K'$ for all $n \in \NN$. Indeed, $K'$ is generated by elements of the form $[k_1,k_2]$ with $k_1,k_2 \in K$, and so, by \cref{equation: generators Kn}, $(K')_n$ is generated by elements of the form $$(1, \dots, 1, [k_1,k_2], 1, \dots, 1)_n$$ with $i = 1, \dots, d^n$. Now, 
\begin{align*}
(1, \dots, 1, [k_1,k_2], 1, \dots, 1)_n =
[(1, \dots, 1, k_1, 1, \dots, 1)_n, \, (1, \dots, 1, k_2, 1, \dots, 1)_n]
\end{align*}
where both elements on the right-hand side are in $K_n$. Since $K_n \leq K$, we conclude that $$(1, \dots, 1, [k_1,k_2], 1, \dots, 1)_n \in K'$$ and the inclusion $(K')_n \leq K'$ follows. 

If $n \geq n_0$, we then have $(K')_{n-n_0} \leq K'$ and therefore $$\pi_{n+1}((K')_{n-n_0}) \leq \pi_{n+1}(K').$$

We are now ready to prove that $\pi_{n+1}(\St_G(n)) \leq \pi_{n+1}(K')$. Let $g \in \pi_{n+1}(\St_G(n))$. By hypothesis, $K \geq \St_G(m)$ and $n_0 \geq m$, so $K \geq \St_G(n_0)$ and by \cref{lemma: St(n+1) = prod St(n)}, the element $g$ can be represented as $$g = (g_1, \dots, g_N)_{n-n_0}$$ with each $g_i \in \pi_{n_0+1}(\St_G(n_0))$. By hypothesis, each $g_i \in \pi_{n_0+1}(K')$, so the element $g \in \pi_{n+1}((K')_{n-n_0}) \leq \pi_{n+1}(K')$ and the result follows.
\end{proof}

To conclude this subsection, we review an equivalence criterion for a regular branch group to be just-infinite.

\begin{Proposition}[{\cite[Proposition 3.5]{BartholdiGrigorchuk2001}}]
Let $T$ be a $d$-regular rooted tree and let $G$ be a subgroup of $ \Aut(T)$ such that $G$ is a regular branch group branching over $K$. Then, $G$ is just-infinite if and only if $K/K'$ is finite. 
\label{proposition: regular branch groups just-infinite equivalence}
\end{Proposition}

\subsection{Congruence subgroup property} \label{section: congruence subgroup property}

Let $T$ be a $d$-regular rooted tree and let $G$ be a subgroup of $\Aut(T)$. Since $\Aut(T)$ is a topological group with $\set{\St(n): n \in \NN}$ as a base of neighborhoods of the identity, then the closure of $G$ in $\Aut(T)$ is isomorphic to $$\overline{G} = \varprojlim G/\St_G(n)$$ by \cite[Corollary 1.1.8]{RibesZalesskii2000}. In the case that $G$ is a branch group, then $$\widetilde{G} := \varprojlim G/\RiSt_G(n)$$ is also a profinite group and since $\RiSt_G(n) \leq \St_G(n)$ we have a well-defined surjective homomorphism $$\psi_1: \widetilde{G} \rightarrow \overline{G}.$$ The kernel of $\psi_1$ is called the \textit{rigid kernel}.

Another standard construction in group theory is the \textit{profinite completion} of a group, defined as the inverse limit $$\widehat{G} := \varprojlim G/N,$$ where the limit runs over all normal subgroups $N$ of finite index in $G$. Since $G$ is branch, the family of rigid stabilizers is a subset of the family of normal subgroups of finite index of $G$, and therefore we have a well-defined surjective homomorphism $$\psi_2: \widehat{G} \rightarrow \widetilde{G}.$$ The kernel of $\psi_2$ is called the \textit{branch kernel}. Composing the last two maps, we obtain an epimorphism $$\psi_3 = \psi_1 \circ \psi_2: \widehat{G} \rightarrow \overline{G}.$$ The kernel of $\psi_3$ is called the \textit{congruence kernel}. 

A group $G$ is said to have the \textit{congruence subgroup property} if the map $\psi_3$ is an isomorphism. 

Using the definition of the inverse limit, it is not hard to see the following equivalences:

\begin{Lemma}
Let $T$ be a $d$-regular rooted tree and $G \leq \Aut(T)$. Then
\begin{enumerate}
\item The rigid kernel is trivial if and only if for any $n \in \NN$, there exists $m \in \NN$ such that $\RiSt_G(n) \geq \St_G(m)$.
\item The branch kernel is trivial if and only if for any $N \lhd_f G$, there exists $n \in \NN$ such that $N \geq \RiSt_G(n)$.
\item The congruence kernel is trivial if and only if for any $N \lhd_f G$, there exists $n \in \NN$ such that $N \geq \St_G(n)$.
\end{enumerate}
\label{lemma: equivalences inclusion profinite kernels}
\end{Lemma}

\section{The group $\IMG(z^2+i)$}
\label{section: The group IMG(z2+i)}

\begin{figure}
\centering
\begin{tikzpicture}[shorten >=1pt,node distance=3cm,on grid,>={Stealth[round]},
    every state/.style={draw=blue!50,very thick,fill=blue!20}, bend angle = 15]

  \node[state] (1)  {$1$};
  \node[state] (a) [right=of 1] {$a$};
  \node[state] (b) [below =of a] {$b$};
  \node[state] (c) [below =of 1] {$c$};

  \path[->] 
  (1) edge [loop above] node  {x/x} ()       
  (a)  edge node [above] {x/$\overline{x}$} (1)
  (b)  edge node [right] {0/0} (a)
  (b)  edge [bend left] node [below] {1/1} (c)
  (c)  edge [bend left] node [above] {0/0} (b)
  (c)  edge node [left] {1/1} (1);
\end{tikzpicture}
\caption{Automaton corresponding to $\IMG(z^2+i)$.}
\label{figure: automaton IMG(z2+i)}
\end{figure}
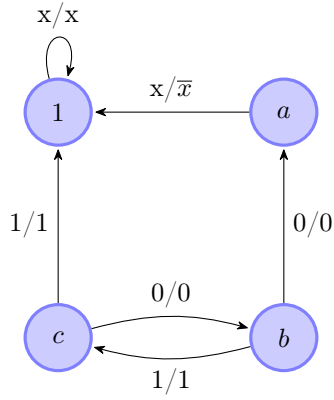

Using the isomorphism $\psi_1$ of \cref{equation: Aut and semidirect product}, we can express the generators of the group $\IMG(z^2+i)$ as 
\begin{equation}
a = (1,1)\sigma, \quad b = (a,c),\text{ and } \quad c = (b,1),
\label{equation: generators IMG}
\end{equation}
where $\sigma$ is the nontrivial permutation in $\Sym(2)$ and $1$ represents the identity element in $\Aut(T)$. This group is also given by the Mealy automaton in \cref{figure: automaton IMG(z2+i)}.

To simplify notation, and since we will only work with this group throughout this article, we refer to it as $G$ from now and on. 

\begin{Lemma}
The elements $a,b,c$ have order $2$. The element $ac$ has order $4$, and the subgroup $\<a,c>$ is isomorphic to $D_4$, the dihedral group with eight elements. 
\label{lemma: order basic elements IMG(z2+i)}
\end{Lemma}

\begin{proof}
The proof of the first statement is proved by induction on the level of the tree. In the first level, $a$ acts like a permutation of order $2$, whereas $b$ and $c$ act trivially. Then, assuming the result holds until level $n$, we have $a^2 = (1,1) = 1$. On the other hand, $b^2 = (a^2,c^2)$ but by the inductive hypothesis $a^2$ and $c^2$ are $1$ until level $n$, so $b^2 = 1$ until level $n+1$. Similarly, $c^2 = (b^2,1)$ and $b^2 = 1$ until level $n$.  

To prove that $ac$ has order $4$, we use the rule with sections presented in \cref{equation: properties sections}:
\begin{align*}
(ac)^4 = ((aca)c)^2 = (b,b)^2 = (b^2,b^2) = (1,1) = 1
\end{align*}

Lastly, we have that $\<a,c> = \<ac, a>$ where $\abs{ac} = 4$, $\abs{a} = 2$ and $$a(ac)a = ca = (ac)^{-1},$$ proving that $\<a,c>$ is the dihedral group with eight elements.
\end{proof}

%\begin{Lemma}
%The group $\IMG(z^2+i)$ is not torsion.
%\end{Lemma}

%\begin{proof}
%Suppose that the element $bca$ has finite order. Then $abc$ and $cab$ have to have the same order. On the other hand, $$(bca)^2 = (abc, cab),$$ which implies that $\abs{bca} = 2 \mcm \set{\abs{abc}, \abs{cab}}$, getting a contradiction. Thus, $bca$ has infinite order.  
%\end{proof}

In \cite[Theorem 2.5]{GrigSavchukSunic2007}, it is proved that a regular branching subgroup of $G$ is $$K := \<[a,b], [b,c]>^G.$$ Moreover, $$G/K \simeq C_2 \times D_4,$$ where the first factor is generated by the class of $b$ and the second factor by the classes of $a$ and $c$. Using \cite[Theorem 9.5]{Radi2025finitetype}, one can in fact prove that $K$ is the maximal regular branching subgroup of $G$.

Let $S_G := \set{a,b,c}$ be the natural generating set of $G$. Define the morphism of free monoids $\phi: \bigoplus S_G \rightarrow \bigoplus S_G$ given by
\begin{align*}
\phi: \left \{ \begin{matrix} a & \mapsto & b, \\
b & \mapsto & c, \\
c & \mapsto & aba.
\end{matrix}\right.
\end{align*}

\begin{Theorem}[{\cite[Theorem 3.1]{GrigSavchukSunic2007}}]
\label{theorem: IMG presentation}
An $L$-presentation for $G$ is
\begin{align*}
G = \langle a,b,c \mid \phi^n(a^2), \phi^n((ac)^4), \phi^n([c,ab]^2), \phi^n([c,bab]^2), \\ \phi^n([c,ababa]^2), \phi^n([c,ababab]^2), \phi^n([c,bababab]^2), n \geq 0 \rangle.
\end{align*}

\end{Theorem}

\section{Generators and presentation of the branching subgroup}
\label{section: generators and presentation of the branching subgroup}

To find a generating set of $K$, we use the Reidemeister algorithm described in \cite[Theorem 2.8 page 91]{CombinatorialGroupTheory}.

First, we select a set of representatives $\mathcal{T}$ of $G/K$:
\begin{align*}
\begin{matrix}
1 & a & ac & aca & acac & cac & ca & c \\
b & ba & bac & baca & bacac & bcac & bca & bc \\
\end{matrix}
\end{align*}

Notice that the first row corresponds to elements in the subgroup $\<a,c> \simeq D_4$ and the second row to the coset $b D_4$ in $G/K$. By choice, each representative has the form $b^\epsilon u$ with $u \in \<a,c>$ and $\epsilon \in \set{0,1}$. We will use this description in subsequent proofs. 

Define the function $\mathscr{R}$ that given a word in $S_G = \set{a,b,c}$, returns the corresponding representative in $\mathcal{T}$. Given $\alpha \in S_G$ and $t \in \mathcal{T}$, define $$s_{t,\alpha} := t\alpha \, \mathscr{R}(t \alpha)^{-1}.$$

Clearly, $s_{t,\alpha}$ belongs to $K$. By \cite[Theorem 2.7 page 89]{CombinatorialGroupTheory}, the set $$\set{s_{t,\alpha}: t \in \mathcal{T}, \alpha \in S_G}$$ generates $K$. 

\begin{Lemma}
If $\alpha \in \set{a,c}$ and $t \in \mathcal{T}$, then $s_{t,\alpha} = 1$. If $u \in \mathcal{T} \cap \<a,c>$, then $s_{bu,b} = (s_{u,b})^{-1}$ and $s_{u,b} = ubu^{-1}b$. 
\label{lemma: stalpha}
\end{Lemma}

\begin{proof}
Write $t = b^\epsilon u$, where $\epsilon \in \set{0,1}$ and $u \in \mathcal{T} \cap \<a,c>$. If $\alpha \in \set{a,c}$, then $t \alpha = b^\epsilon u'$, with $u'$ another word in $\<a,c>$. Therefore $\mathscr{R}(t \alpha) = b^\epsilon u'$ and $$s_{t, \alpha} = b^\epsilon u' (u')^{-1} b^{\epsilon} = 1.$$

Notice that $\mathscr{R}(ub) = bu$ and $\mathscr{R}(bub) = u$, so $$s_{u,b} = ub u^{-1}b$$ and
\begin{equation*}
s_{bu,b} = bub u^{-1} = (ubu^{-1}b)^{-1} = (s_{u,b})^{-1}. \qedhere
\end{equation*}
\end{proof}

We can also find the relations in $K$ in terms of the generators $s_{t,\alpha}$. Following \cite{CombinatorialGroupTheory}, we define a function $\tau$ that transforms words in $S_G$ into words in the set $\set{s_{t,\alpha}: t \in T, \alpha \in S_G}$ as follows: if $\omega = \alpha_1 \dots \alpha_n$ is a word in $S_G$, then define $\tau(\omega) := s_{t_1,\alpha_1} \dots s_{t_n, \alpha_n}$ where $t_i = \mathscr{R}(\alpha_1 \dots \alpha_{i-1}) \in \mathcal{T}$.

For example, consider the word $abc$. Then $$\tau(abc) = s_{1,a}s_{a,b}s_{ba,c}.$$

\begin{Theorem}[{\cite[Theorem 2.8]{CombinatorialGroupTheory}}]
A presentation for $K$ is given by $$\set{s_{t,\alpha}: t \in \mathcal{T}, \alpha \in S_G}$$ as generators and the relations are $$\tau(s_{t, \alpha})s_{t, \alpha}^{-1} \text{ and } \tau(tRt^{-1}),$$ where $\alpha \in S_G$, $t \in \mathcal{T}$ and $R$ is a relation of the presentation of $G$.
\label{theorem: presentation subgroup algorithm}
\end{Theorem}

Define the elements 
\begin{itemize}
\item $x = [a,b]$
\item $y = [b,c]$
\item $z = aya$
\item $t = cxcy^{-1}x^{-1}$
\item $w = ata$
\end{itemize}

Using the description by sections of $a,b$ and $c$, it is not hard to see that
\begin{itemize}
\item $y = (x,1)$
\item $z = (1,x)$
\item $t = (y,1)$
\item $w = (1,y)$
\end{itemize}

\begin{Lemma}
The subgroup $K = \<x,y,z,t,w>$.
\label{lemma: K is xyztw}
\end{Lemma}

\begin{proof}
Recall that $K = \<[a,b],[b,c]>^G = \<x,y>^G$. Then, by the definitions of $x,y,z,w$ and $t$, it is clear that $\<x,y,z,w,t> \subseteq K$. To prove the other inclusion, we use \cref{lemma: stalpha} and the fact that $$K = \<s_{t, \alpha}: t \in T, \alpha \in S_G> = \<s_{u,b}: \text{$u \in \mathcal{T} \cap \<a,c>$}>.$$ We examine the eight possible values of $u$ and prove that in each case, we can write $s_{u,b}$ in terms of $x,y,z,w$ and $t$. 

If $u = 1$, then $s_{1,b} = 1$. 

If $u = a$, then $s_{a,b} = abab = x$. 

If $u = ac$, then $s_{ac,b} = acbcab = z^{-1}x$.

If $u = aca$, then $s_{aca,b} = acabacab = w$.

If $u = acac$, then $s_{acac,b} = acacbcacab = y^{-1}w$.

If $u = cac$, then $s_{cac,b} = cacbcacb = z^{-1}tx$.

If $u = ca$, then $s_{ca,b} = cabacb = tx$.

If $u = c$, then $s_{c,b} = cbcb = y^{-1}$.
\end{proof}

\section{The exponent function in $K$}
\label{section: The exponent function in K}

Define $S_K = \set{x,y,z,t,w}$ and the function $e: \bigoplus S_K \rightarrow \ZZ/2\ZZ$ where given a word with letters in $S_K$, it returns the parity of the length of the word. Clearly, the function $e$ is a morphism from the free monoid $\bigoplus S_K$ to $\ZZ/2\ZZ$. In this section, we aim to prove that $e$ descends to a well-defined function $e: K \rightarrow \ZZ/2\ZZ$. To achieve this, we must show that for every relation $\mathfrak{R}$ in $S_K$, we have $e(\mathfrak{R}) = 0$.

The simplest approach would be to imitate the proof for the Basilica group presented in \cite[Lemma 1 page 7]{GrigZuk2002}. However, there is a big difference between that case and ours. In the Basilica group, if $g = (g_0,g_1)$, then the lengths of $g_0$ and $g_1$ in terms of the generators are strictly smaller than the length of $g$. In our case, this property does not hold, as for example $xy^{-1}x^{-1} = (txy,1)$, where the length remains the same or $xy^{-2}x^{-1} = (txytxy,1)$, where the length is actually bigger. For this reason, we adopt this more sophisticated strategy, which involves finding all relations $\mathfrak{R}$ in $K$ and verifying that $e(\mathfrak{R}) = 0$. Thus, we need to show that $$e(\tau(s_{t, \alpha})s_{t, \alpha}^{-1}) = e(\tau(tRt^{-1})) = 0 \pmod{2}$$ for all $\alpha \in S_G$, $t \in \mathcal{T}$ and $R$ a relation in $G$. 

Define also the function $\ell: \bigoplus S_G \rightarrow \ZZ/2\ZZ$ where given a word in $S_G$, it returns the parity of the number of letters $a$ and $c$ in the word and $\beta: \bigoplus S_G \rightarrow \ZZ/2\ZZ$, that returns the parity of the number of letters $b$ in the word. For example, we have $\ell(abacac) = 5 \equiv 1 \pmod{2}$ and $\beta(abacac) \equiv 1 \pmod{2}$.

\begin{Lemma}
The functions $\ell$ and $\beta$ descend to well-defined functions from the quotient $G/K$ to $\ZZ/2\ZZ$. Moreover, if $\omega$ is a word in $S_G$, then we have $\ell(\omega^n) = n \ell(\omega)$ and $\beta(\omega^n) = n \beta(\omega)$.
\label{lemma: ell beta descend to G/K}
\end{Lemma}

\begin{proof}
The first claim follows from the fact that $K = \<[a,b],[b,c]>^G$ and the generators satisfy $\beta([a,b]) = \beta([b,c]) = \ell([a,b]) = \ell([b,c]) = 0$. The second part is straightforward. 
\end{proof}

\begin{Lemma}
For all $t \in \mathcal{T}$, we have $e(s_{t,b}) = \ell(t)$.
\label{lemma: relation e and l}
\end{Lemma}

\begin{proof}
Consider first the case when $t = u$ with $u \in \mathcal{T} \cap \<a,c>$. In this case, the result follows directly from the explicit expressions for the elements $s_{u,b}$ given in \cref{lemma: K is xyztw}. 

Now consider the case when $t = bu$ with $u \in \mathcal{T} \cap \<a,c>$. Using the property $s_{bu,b} = (s_{u,b})^{-1}$, we have $$e(s_{bu,b}) = e((s_{u,b})^{-1}) = e(s_{u,b}) = \ell(u) = \ell(bu),$$ completing the proof.
\end{proof}

\begin{Lemma}
\label{lemma: formula e(tau(tRt-1))}
Let $t \in \mathcal{T}$ and $R$ be a word in $S_G$. Then $$e(\tau(tRt^{-1})) = \beta(R)\ell(t) + e(\tau(R)) + \beta(t)\ell(R).$$   
\end{Lemma}

\begin{proof}
Write $t = \beta_1\dots\beta_m$ and $R = \alpha_1\dots\alpha_n$. Then, $$tRt^{-1} = \beta_1\dots\beta_m \alpha_1\dots\alpha_n \beta_m\dots\beta_1.$$ Applying the function $e$ to $\tau(tRt^{-1})$ we obtain 
\begin{align*}
e(\tau(tRt^{-1})) = \sum_{i = 1}^m e(s_{\mathscr{R}(\beta_1\dots\beta_{i-1}),\beta_i}) + \sum_{j = 1}^n e(s_{\mathscr{R}(t\alpha_1\dots\alpha_{j-1}),\alpha_j}) \\ + \sum_{k = 1}^m e(s_{\mathscr{R}(tR\beta_m\dots\beta_{k+1}),\beta_k}).
\end{align*}

By the choice of the elements in $\mathcal{T}$, only $\beta_1$ can be $b$, and by \cref{lemma: stalpha},  $s_{*,\alpha} = 1$ for $\alpha \in \set{a,c}$. This simplifies the previous formula to
\begin{align*}
e(\tau(tRt^{-1})) = e(s_{1,\beta_1}) + \sum_{j: \alpha_j = b} e(s_{\mathscr{R}(t\alpha_1\dots\alpha_{j-1}),b}) + e(s_{\mathscr{R}(tR\beta_m\dots\beta_2),\beta_1}).
\end{align*}

We have $s_{1,\alpha} = 1$ for all $\alpha \in S_G$, $e(s_{\mathscr{R}(\omega),b}) = \ell(\mathscr{R}(\omega))$ by \cref{lemma: relation e and l} and 
$\ell(\mathscr{R}(\omega)) = \ell(\omega)$ by \cref{lemma: ell beta descend to G/K}. Thus, 
\begin{align*}
e(\tau(tRt^{-1})) = \sum_{j: \alpha_j = b} \ell(t\alpha_1\dots\alpha_{j-1}) + e(s_{\mathscr{R}(tR\beta_m\dots\beta_2),\beta_1}).
\end{align*}

Moreover, by the linearity of $\ell$,
\begin{align*}
e(\tau(tRt^{-1})) = \sum_{j: \alpha_j = b} \left( \ell(t) + \ell(\alpha_1\dots\alpha_{j-1}) \right) + e(s_{\mathscr{R}(tR\beta_m\dots\beta_2),\beta_1}) \\
= \ell(t) \beta(R) + \sum_{j: \alpha_j = b} \ell(\alpha_1\dots\alpha_{j-1}) + e(s_{\mathscr{R}(tR\beta_m\dots\beta_2),\beta_1}).
\end{align*}

On the other hand, 
\begin{align*}
e(\tau(R)) = \sum_{j = 1}^n e(s_{\mathscr{R}(\alpha_1\dots\alpha_{j-1}),\alpha_j}) = \sum_{j: \alpha_j = b} \ell(\alpha_1\dots\alpha_{j-1}),
\end{align*}
so combining the last two equations, 
\begin{align*}
e(\tau(tRt^{-1})) = \ell(t) \beta(R) + e(\tau(R)) + e(s_{\mathscr{R}(tR\beta_m\dots\beta_2),\beta_1}).
\end{align*}

Finally, we evaluate the last term. If $\beta_1 \neq b$, then $e(s_{\mathscr{R}(tR\beta_m\dots\beta_2),\beta_1}) = 0$. If $\beta_1 = b$, 
\begin{align*}
e(s_{\mathscr{R}(tR\beta_m\dots\beta_2),b}) = \ell(tR\beta_m\dots\beta_2) = \ell(t) + \ell(R) + \ell(\beta_m\dots\beta_2) \\ = \ell(t) + \ell(R) + \ell(t^{-1}) = \ell(R),   
\end{align*}
where in the last two steps it was used that $\beta_1 = b$, so $$\ell(\beta_m\dots\beta_2) = \ell(\beta_m\dots\beta_1) = \ell(t^{-1}) = \ell(t).$$ Therefore, combining both cases, we conclude
\begin{equation*}
e(s_{\mathscr{R}(tR\beta_m\dots\beta_2),\beta_1}) = \beta(t) \ell(R). \qedhere
\end{equation*}
\end{proof}

\begin{Lemma}
\label{lemma: IMG e(tau(tRt-1)) for squared words}
Let $R$ be a word in $S_G$ and $t \in \mathcal{T}$. Suppose there exists another word $R'$ in $S_G$ such that $R = R'^2$. Then, $$e(\tau(tRt^{-1})) = e(\tau(R)) = \ell(R') \beta(R').$$
\end{Lemma}

\begin{proof}
First, by \cref{lemma: ell beta descend to G/K} and \cref{lemma: formula e(tau(tRt-1))}, we have 
\begin{align*}
e(\tau(tRt^{-1})) = \beta(R)\ell(t) + e(\tau(R)) + \beta(t)\ell(R) = \\ \beta(R'^2)\ell(t) + e(\tau(R)) + \beta(t)\ell(R'^2) = 2\beta(R')\ell(t) + e(\tau(R)) + 2\beta(t)\ell(R') \\ = e(\tau(R)).
\end{align*}

\bigskip

Next, write $$R' = \alpha_{11}\dots\alpha_{1m_1}b^{r_1}\alpha_{21}\dots\alpha_{2m_2}b^{r_2}\dots\alpha_{n1}\dots\alpha_{nm_n}b^{r_n},$$ where $m_i \neq 0$ for $i > 1$, $\alpha_{ij} \in \set{a,c}$ and $r_k \neq 0$ for $k \neq n$. The number $e(\tau(R))$ is obtained by adding up the values of $\ell(\omega_k(R))$'s, where $\omega_k(R)$ is the segment of $R$ preceding a $b$. Since $R = R'^2$, this yields the formula:
\begin{align*}
e(\tau(R)) = r_1m_1 + r_2(m_1 + m_2) + \dots + r_n(m_1+\dots+m_n) + \\ r_1(m_1+\dots+m_n+m_1) + \dots + r_n(m_1 + \dots +m_n + m_1 + \dots + m_n).
\end{align*}

Now, $m_1 + \dots + m_n = \ell(R')$ and $r_1 + \dots + r_n = \beta(R')$, so the expression simplifies to
\begin{align*}
e(\tau(R)) = 2 \sum_{k = 1}^n r_k \left(\sum_{i = 1}^k m_i \right) + (r_1 + \dots + r_n) \ell(R') = \ell(R') \beta(R')
\end{align*}
as we wanted.
\end{proof}

Recall the substitution $\phi$ presented in \cref{section: The group IMG(z2+i)}, that defines the presentation of $G$. This function $\phi: \bigoplus S_G \rightarrow \bigoplus S_G$ is given by
\begin{align*}
\phi: \left \{ \begin{matrix} a & \mapsto & b \\
b & \mapsto & c \\
c & \mapsto & aba.
\end{matrix}\right.
\end{align*}

\begin{Lemma}
Let $R$ be a word in $S_G$. Then $\ell(\phi(R)) = \beta(R)$ and $\beta(\phi(R)) = \ell(R)$.
\label{lemma: IMG ell beta and phi}
\end{Lemma}

\begin{proof}
The result is immediate from the following three facts: each occurrence of $a$ is replaced by $b$. Each occurrence of $b$ is replaced by $c$. Each occurrence of $c$ is replaced by two $a$'s (which do not change $\ell$) and one $b$.
\end{proof}

We are finally ready to prove the main lemma of this section:

\begin{Lemma}
The function $e: \bigoplus S_K \rightarrow \ZZ/2\ZZ$ descends to a well-defined function on $K$. 
\label{lemma: e descends to a function in K}
\end{Lemma}

\begin{proof}
As it was previously mentioned, we need to show that if $\mathfrak{R}$ is a relation in $K$, then $e(\mathfrak{R}) = 0$. By \cref{theorem: presentation subgroup algorithm}, the relations in $K$ are of the form $$\tau(s_{t, \alpha})s_{t, \alpha}^{-1} \text{ and } \tau(tRt^{-1}),$$ where $\alpha \in S_G$, $t \in \mathcal{T}$ and $R$ is a relation of the presentation of $G$.

For the first family of relations and by \cref{lemma: stalpha}, we only need to check the case $s_{t,\alpha}$ for $\alpha = b$. Write $t = b^\epsilon u$ with $\epsilon \in \set{0,1}$ and $u \in \mathcal{T} \cap \<a,c>$. 

If $\epsilon = 0$, then $s_{u,b} = ubu^{-1}b$ and so  $$\tau(s_{u,b}) = s_{\mathscr{R}(u),b} s_{\mathscr{R}(ubu^{-1}),b}$$ since $s_{*,\alpha} = 1$ for all $\alpha \in \set{a,c}$ by \cref{lemma: stalpha}. Thus,
\begin{align*}
e(\tau(s_{u,b}) s_{u,b}^{-1}) = e(s_{\mathscr{R}(u),b}) + e(s_{\mathscr{R}(ubu^{-1}),b}) + e(s_{u,b}) = \\ \ell(\mathscr{R}(u)) + \ell(\mathscr{R}(ubu^{-1})) + \ell(u) = \ell(u) + \ell(ubu^{-1}) + \ell(u) = 0.  
\end{align*}
Similarly, if $\epsilon = 1$, then $s_{bu,b} = bubu^{-1}$ and so  $\tau(s_{bu,b}) = s_{1,b} s_{\mathscr{R}(bu),b}$. Thus,
\begin{align*}
e(\tau(s_{bu,b}) s_{bu,b}^{-1}) = e(s_{\mathscr{R}(bu),b}) + e(s_{bu,b}) = \\ \ell(\mathscr{R}(bu)) + \ell(bu) = \ell(bu) + \ell(bu) = 0.  
\end{align*}

We now prove the same for the relations of the form $\tau(tRt^{-1})$, where $t \in \mathcal{T}$ and $R$ is a relation of $G$. By \cref{theorem: IMG presentation} all relations have the form $\phi^n(R'^2)$ for some word $R'$. Moreover, since $\phi$ is a morphism of monoids, the relations can be written in the form $\phi^n(R')^2$. Applying \cref{lemma: IMG e(tau(tRt-1)) for squared words}, we obtain
\begin{align*}
e(\tau(t\phi^n(R')^2t^{-1})) = e(\tau(\phi^n(R')^2)) = \ell(\phi^n(R')) \beta(\phi^n(R')).
\end{align*}
By \cref{lemma: IMG ell beta and phi}, the last quantity is invariant under $\phi$, so 
\begin{align*}
e(\tau(t\phi^n(R')^2t^{-1})) = \ell(R') \beta(R'),
\end{align*}
reducing the problem to a finite verification. The words $R'$ that appear in the presentation of $G$ are: $$a, \, (ac)^2, \, [c,ab], \, [c,bab], \, [c,ababa], \, [c,ababab] \text{ and } [c,bababab].$$ Each of these words has an even number of ocurrences of $b$, meaning that $\beta(R') = 0$ for all of them, proving the lemma.
\end{proof}

\section{Abelianization of $K$}
\label{section: abelianization of K}

In this section, we will prove \cref{theorem: K/K' is C45}, namely, that $K/K' \simeq C_4^5$.

\begin{Lemma}
We have $[K:K'] \leq 4^5$. In particular $G$ is just-infinite. 
\label{lemma: [K:K'] leq 64}
\end{Lemma}

\begin{proof}
As it was proved in \cref{lemma: K is xyztw}, the group $K$ is generated by five elements and it is not hard to see that all of them have order $4$. Therefore, $K/K'$ can be embedded in $C_4^5$. Since $K/K'$ is finite, the group $G$ is just-infinite by \cref{proposition: regular branch groups just-infinite equivalence}.
\end{proof}

\begin{Lemma}
The following conjugation relations hold: 
\begin{align*}
\begin{matrix}
axa = x^{-1} & bxb = x^{-1} & cxc = txy \\
aya = z & byb = y^{-1} & cyc = y^{-1} \\
aza = y & bzb = x^{-1}z^{-1}x & czc = z \\
ata = w & btb = x^{-1}t^{-1}x & ctc = t^{-1} \\
awa = t & bwb = w^{-1} & cwc = w \\
\end{matrix}
\end{align*}
In particular, for any $\alpha \in S_G$ and $\beta \in S_K$, we have $e(\alpha \beta \alpha) = 1$.
\label{lemma: relations conjugations K IMG}
\end{Lemma}

\begin{proof}
We label the relations as in a $5 \times 3$-matrix format for clarity. Observe that $axa = aababa = baba = x^{-1}$. The relations 11, 12, 22 and 23 follow in the same way.

The relations 21, 31, 41, 51 and 13 follow directly from the definition of the elements.

The relations 52, 33, 43, and 53 are obtained by applying the previous relations on the first level. For example $$bwb = (a,c)(1,y)(a,c) = (1, cyc) = (1, y^{-1}) = w^{-1}.$$

For 32, $$x^{-1}z^{-1}x = (baba)(ay^{-1}a)(abab) = baby^{-1}bab = bayab = bzb.$$

For 42, $$btb = (bcxcb)(by^{-1}b)(bx^{-1}b) = (bcxcb)yx.$$ Then, $$bcxcb = bcababcb = (bcabac)(cbcb) = x^{-1}t^{-1}y^{-1},$$ and hence $btb = x^{-1}t^{-1}x$.
\end{proof}

\begin{Lemma}
We have $K \geq \St_G(3) \geq K'$.
\label{lemma: K contains St(3)}
\end{Lemma}

\begin{proof}
We first prove that $\St_G(3) \geq K'$. By \cref{lemma: commutator generators}, we have $$K' = \<[\alpha,\beta]: \alpha, \beta \in \set{x,y,z,w,t}>^G.$$ As $w,t \in \St_G(3)$ and $\St_G(3) \lhd G$, it is enough to observe that the commutators $[x,y], [x,z], [y,z] \in \St_G(3)$.

We offer two proofs of $K \geq \St_G(3)$. For the first one, we apply \cref{lemma: map defined in abelianizations} with $(A,B, \pi) = (G, K, \pi_3)$. Using GAP \cite{GAP}, we find that $$[G:K] = [\pi_3(G):\pi_3(K)] = 16$$ and so $\varphi$ is an isomorphism, which implies that $\ker(\pi_3) = \St_G(3) \leq K$.

The second proof is to observe that $K \geq \St_G(3)$ is equivalent to satisfying that for all $\gamma \in \mathcal{T} \setminus \set{1}$ then $$\gamma x^{\alpha_1} y^{\alpha_2} z^{\alpha_3} \notin \St_G(3)$$ for all $\alpha_1 \in \set{0, \dots, 3}$ and $\alpha_2, \alpha_3 \in \set{0, 1}$.

Indeed, if $K \ngeq \St_G(3)$, then $\St_G(3) \cap K^c \neq \emptyset$ and so there exists $\gamma \in \mathcal{T} \setminus \set{1}$ such that $\gamma K \cap \St_G(3) \neq \emptyset$. Let $k \in K$ such that $\gamma k \in \St_G(3)$. As $K$ is generated by $\set{x,y,z,t,w}$ and  $\<K', t,w, y^2, z^2> \subseteq K \cap \St_G(3)$, we may multiply by an element $k' \in \<K', t,w,y^2,z^2>$ such that $k k'$ is of the form $x^{\alpha_1} y^{\alpha_2} z^{\alpha_3}$ and $\gamma x^{\alpha_1} y^{\alpha_2} z^{\alpha_3} \in \St_G(3)$. 

Suppose there exist $\gamma \in \mathcal{T} \setminus \set{1}$, $\alpha_1 \in \set{0, \dots, 3}$, $\alpha_2, \alpha_3 \in \set{0,1}$ such that $\gamma x^{\alpha_1} y^{\alpha_2} z^{\alpha_3} \in \St_G(3)$. As $x,y,z \in \St_G(1)$, then $\gamma$ must be in $\St_G(1)$. Write $\gamma = b^\epsilon u$ with $\epsilon \in \set{0,1}$ and $u \in \set{1, aca, acac, c}$. As the elements $u,x^2,y,z \in \St_G(2)$ then $$b^\epsilon x^\delta = (a^\epsilon (ca)^\delta, c^\epsilon (ac)^\delta)$$ must be in $\St_G(2)$, where  $\epsilon, \delta \in \set{0,1}$. This implies that $\epsilon = \delta = 0$ and so $$\gamma x^{\alpha_1} y^{\alpha_2} z^{\alpha_3} = u (x^2)^{\beta_1} y^{\alpha_2} z^{\alpha_3}$$ with $u \in \set{aca,acac,c}$ and $\beta_1, \alpha_2, \alpha_3 \in \set{0,1}$. 

If $u = c$, then $$ u (x^2)^{\beta_1} y^{\alpha_2} z^{\alpha_3} = (ab^{\beta_1} (ca)^{\alpha_2}, c b^{\beta_1} (ac)^{\alpha_2}, b^{\beta_1} (ca)^{\alpha_3}, b^{\beta_1} (ac)^{\alpha_3})_2$$ must be in $\St_G(3)$, or equivalently, every coordinate must be in $\St_G(1)$. By the first coordinate, $\alpha_2 = 1$ and by the second coordinate $\alpha_2 = 0$, leading to a contradiction. The other cases of $u$ follow by the same argument. This concludes the proof of the inclusion $K \geq \St_G(3)$. 
\end{proof}

Recall $\gamma_3(G) := [G',G]$. Since $G$ is generated by $a$, $b$ and $c$, then $$\gamma_3(G) = \<[\alpha, \beta, \gamma]: \alpha, \beta, \gamma \in \set{a,b,c}>^G.$$

Moreover, we only need the following commutators:

\begin{itemize}
\item $[a,b,a] = x^2$,
\item $[a,b,b] = x^2$,
\item $[a,b,c] = xy^{-1}x^{-1}t^{-1}$,
\item $[a,c,a] = 1$,
\item $[a,c,b] = y^{-1}w$,
\item $[a,c,c] = 1$,
\item $[b,c,a] = yz^{-1}$,
\item $[b,c,b] = y^2$,
\item $[b,c,c] = y^2$,
\end{itemize}

As $w[a,c,b]^{-1}w^{-1} = yw^{-1}$, we conclude that $$\gamma_3(G) = \<x^2,y^2,txyx^{-1},yw^{-1}, yz^{-1}>^G.$$

\begin{Lemma}
We have $K' \geq (\gamma_3(G))_2$.
\label{lemma: K' contains gamma3(G)2}
\end{Lemma}

\begin{proof}
Using the relations in \cref{lemma: relations conjugations K IMG} we have the following identities:
\begin{align*}
[y,(ca)^2] = y^2, \quad [x,ca] = xtxy, \quad [y,ca] = yz^{-1}, \quad [t,ca] = tw^{-1}.
\end{align*}

Then, the element $[t,x^2] \in K'$ and $$[t,x^2] = ([y,(ca)^2],1) = (y^2,1) = (x^2,1,1,1)_2.$$

Next, as $G$ is regular branch over $K$, the elements $(t,1), (x^2,1), (t,1,1,1)_2 \in K$. Therefore  $$[(t,1),(x^2,1)] = ([t,x^2],1) = (y^2,1,1,1)_2 \in K',$$  $$[(t,1),y] = ([t,x],1) = ([y,ca],1,1,1)_2 = (yz^{-1},1,1,1)_2 \in K'$$ and $$[(t,1,1,1)_2,y] = ([t,ca],1,1,1)_2 = (tw^{-1},1,1,1)_2 \in K'.$$ 
As $(K')_2 \leq K'$ and $(txyx^{-1},1,1,1)_2 \in K'$, we can permute the elements to conclude that $(ty,1,1,1)_2 \in K'$. Taking the inverse and multiplying by $(y^2,1,1,1)_2$ and $(tw^{-1},1,1,1)_2$ gives $$(y^2 y^{-1}t^{-1} tw^{-1},1,1,1)_2 = (yw^{-1},1,1,1)_2 \in K'.$$
Also $$[t,y] = ([y,x],1) = ([x,ca],1,1,1)_2 = (xtxy,1,1,1)_2 \in K'.$$
Conjugating the above by $(x^{-1},1,1,1)_2 \in K$ yields $(txyx,1,1,1)_2 \in K'$ and since $(x^{-2},1,1,1)_2 \in K'$ we conclude $$(txyx^{-1},1,1,1)_2 \in K'.$$

As $G$ is fractal, this implies that $K' \geq (\gamma_3(G))_2$.
\end{proof}

\begin{Lemma}
We have $\pi_4(\gamma_3(G)) \geq \pi_4(\St_G(3))$ and $\pi_6(K') \geq \pi_6(\St_G(5))$.
\label{lemma: pi4(gamma3) contains pi4(stG(3))}
\end{Lemma}

\begin{proof}
We first observe that $\pi_4(y) \in \pi_4(\gamma_3(G))$. Indeed, one can compute that $$txyx^{-1}y^{-1} = (b, aba, b,b, 1, 1,1, 1)_3$$ and as $b, aba \in \St_G(1)$, then $$\pi_4(txyx^{-1}) = \pi_4(y) \in \pi_4(\gamma_3(G)).$$ 

This simplifies the normal generators of $\pi_4(\gamma_3(G))$ to $$\pi_4(\gamma_3(G)) = \<\pi_4(x^2), \pi_4(y), \pi_4(z), \pi_4(t), \pi_4(w)>^{\pi_4(G)}.$$

To prove $\pi_4(\gamma_3(G)) \geq \pi_4(\St_G(3))$, we apply the same trick done in \cref{lemma: K contains St(3)}. As $K \geq \St_G(3)$, we can use a set of representatives of non-trivial cosets of $K/\gamma_3(G)$ and a set of representatives of cosets of $\gamma_3(G)/(\gamma_3(G) \cap \St_G(3))$. Since we have the inclusion $\<K', t,w, x^2, y^2, z^2> \leq \St_G(3)$, it is enough to prove that $$x y^{\alpha_1} z^{\alpha_2} \notin \St_G(3)$$ for $\alpha_1, \alpha_2 \in \set{0,1}$. As $y, z \in \St_G(2)$ but $x \in \St_G(1) \setminus \St_G(2)$, this is clearly true.  

Finally, combining the previous result with \cref{lemma: K' contains gamma3(G)2} we conclude 
\begin{equation*}
\pi_6(K') \geq \pi_6(\gamma_3(G)_2) = (\pi_4(\gamma_3(G)))_2 \geq (\pi_4(\St_G(3)))_2 = \pi_6(\St_G(5)). \qedhere
\end{equation*}
\end{proof}

\cref{lemma: pi4(gamma3) contains pi4(stG(3))} cannot be generalized to $\gamma_3(G) \geq \St_G(3)$ because $y \notin \gamma_3(G)$. This is basically  because elements in $\gamma_3(G)$ can be written with an even number of generators of $S_K$ by \cref{lemma: e descends to a function in K}. This is one of the facts that explain why $\IMG(z^2+i)$ does not have the congruence subgroup property, since otherwise $K'$ would contain $\St_G(5)$ and by \cite[Proposition 3.8]{BartholdiGrigorchuk2001} then $\IMG(z^2+i)$ would have the congruence subgroup property.

\begin{Lemma}
For all $n \geq 5$ we have $\pi_n(K)^{\ab} \simeq  C_4^3$. The group is generated by the classes of $[\pi_n(x)],[\pi_n(y)]$ and $[\pi_n(z)]$ and we have the relations $[\pi_n(t)] = [\pi_n(y)]^2$ and $[\pi_n(w)] = [\pi_n(z)]^2$. 
\label{lemma: IMG pin(K)ab}
\end{Lemma}

\begin{proof}
By \cref{lemma: K contains St(3)}, we have $K \geq \St_G(3)$. By \cref{lemma: pi4(gamma3) contains pi4(stG(3))}, we also have the inclusion $\pi_6(K') \geq \pi_6(\St_G(5))$. So, by \cref{lemma: stabilization pin(K)ab}, it follows that $$\pi_n(K)^{\ab} \simeq \pi_5(K)^{\ab}$$ for all $n \geq 5$.
Since $K$ is generated by $x,y,z,w$ and $t$, then the abelianization $\pi_n(K)^{\ab}$ is generated by their respective classes. The only relations appearing are $[\pi_5(t)] = [\pi_5(y)]^2$ and $[\pi_5(w)] = [\pi_5(z)]^2$, and therefore $$\pi_5(K)^{ab} \simeq C_4^3,$$ proving the result. 
\end{proof}

Define the subgroup $$\mathcal{K} := \set{(g_0,g_1) \in K_1: e(g_0) = e(g_1) \equiv 0 \pmod{2}}.$$

\begin{Proposition}
The commutator subgroup of $K$ satisfies $K' \leq \mathcal{K}$.
\label{proposition: IMG K' contain in even parity}
\end{Proposition}

\begin{proof}
By \cref{lemma: commutator generators}, we know that $K' = \<[s,s']: s,s' \in S_K>^K$. Moreover, since $K \lhd G$, we may replace the exponent $K$ by any subgroup containing it. For convenience, we take $\St_G(1)$.

We begin by verifying that each commutator $[s,s']$ satisfies the expected parity condition. Since $[s',s] = [s,s']^{-1}$ and $[s,s] = 1$, we only need to check ten commutators:

\begin{itemize}
\item $[x,y] = ([ca,x], 1) = ((caxac)x^{-1},1)$. Since $K$ is normal, $caxac \in K$, so $[ca,x] \in K$. Moreover, by \cref{lemma: relations conjugations K IMG} we have $e(caxac) = e(x) \equiv 1 \pmod{2}$, so $e([ca,x]) \equiv 0 \pmod{2}$. 

\item $[x,z] = (1,[ac,x])$ and the verification follows from the same argument as the previous case. 

\item $[x,t] = ([ca,y],1)$ and the same argument applies.

\item $[x,w] = (1,[ac,y])$ and the same argument applies.

\item $[y,z] = 1$.

\item $[y,t] = ([x,y],1)$, that clearly has an even number of elements of $S_K$. 

\item $[y,w] = 1$.

\item $[z,t] = 1$.

\item $[z,w] = (1, [x,y])$, that clearly has an even number of elements of $S_K$. 

\item $[z,t] = 1$.
\end{itemize}

Now, if $h = (h_0, h_1) \in \St_G(1)$ and $s,s' \in S_K$ with $[s,s'] = (g_0,g_1)$, we have $$h[s,s']h^{-1} = (h_0g_0h_0^{-1}, h_1g_1h_1^{-1}).$$ Since $g_0,g_1 \in K$ and $K \lhd G$, it follows that $h_0g_0h_0^{-1}, h_1g_1h_1^{-1} \in K$. Moreover, again by \cref{lemma: relations conjugations K IMG}, we have $e(h_ig_ih_i^{-1}) = e(g_i) \equiv 0 \pmod{2}$ for $i = 0,1$, proving the result.
\end{proof}

We are now ready to prove \cref{theorem: K/K' is C45}.

\begin{Theorem}
The abelianization of $K$ satisfies $K/K' \simeq C_4^5$.
\end{Theorem}

\begin{proof}
Applying \cref{lemma: map defined in abelianizations} with $(A,B,\pi) = (K, K', \pi_5)$, we obtain a surjective homomorphism $$\varphi: K/K' \rightarrow \pi_5(K)/\pi_5(K').$$ Thus, the only possible additional relations in $K/K'$ must arise from those in $\pi_5(K)/\pi_5(K')$. By \cref{lemma: IMG pin(K)ab}, these potential relations can only be $ty^{-2}$ or $wz^{-2}$ in $K'$. However, neither relation holds because $$ty^{-2} = (yx^{-2},1) \text{ and } wz^{-2} = (1,yx^{-2}),$$ and the parity of the coordinates is not even, meaning that $ty^{-2}$ and $wz^{-2}$ are not in $\mathcal{K}$. Since $\mathcal{K}$ contains $K'$, by \cref{proposition: IMG K' contain in even parity}, it follows that $ty^{-2}$ and $wz^{-2}$ are not in $K'$.
\end{proof}

As a corollary of \cref{theorem: K/K' is C45}, we complete the proof of \cref{theorem: IMG joo and no CSP}.

\begin{Theorem}
The group $G$ does not have the congruence subgroup property.
\end{Theorem}

\begin{proof}
By \cref{lemma: [K:K'] leq 64} and the fact that $G$ is regular branch over $K$, we know that $K'$ has finite index in $G$. So, in order to prove that $G$ does not have the congruence subgroup property, it suffices to show that $K'$ does not contain $\St_G(n)$ for any $n \in \NN$ (see \cref{lemma: equivalences inclusion profinite kernels}). Since the stabilizers subgroups are all nested, we may assume without loss of generality that $n \geq 5$. 

Suppose for a contradiction that $K'$ contains $\St_G(n)$ for some $n \geq 5$. Applying \cref{lemma: map defined in abelianizations} with $(A,B,\pi) = (K, K', \pi_n)$, we obtain that the natural surjective map $\varphi: K^{\ab} \rightarrow \pi_n(K)^{\ab}$ must be an isomorphism. However, from \cref{lemma: IMG pin(K)ab}, we know that $\pi_n(K)^{ab} \simeq C_4^3$ for all $n \geq 5$ and that $K/K' \simeq C_4^5$ by \cref{theorem: K/K' is C45}.
\end{proof}

To conclude this section, we provide additional information about the subgroup $\mathcal{K}$.

\begin{Proposition}
The subgroup $\mathcal{K}$ is normal in $G$. The quotient $K/\mathcal{K}$ is isomorphic to $C_4 \times C_2 \times C_2$ and $[\mathcal{K}: K'] = 64$.
\end{Proposition}

\begin{proof}
To prove that $\mathcal{K} \lhd G$, we need to check that for any $\alpha \in S_G$ and any $(g_0, g_1) \in \mathcal{K}$, we have $\alpha (g_0, g_1) \alpha \in \mathcal{K}$.

If $\alpha = a$, then $a (g_0, g_1) a = (g_1, g_0) \in \mathcal{K}$ since the condition is imposed over both coordinates.

If $\alpha = b$, then $b (g_0, g_1) b = (a g_0 a, c g_1 c)$ and by \cref{lemma: relations conjugations K IMG}, $e(a g_0 a) = e(g_0)$ and $e(c g_1 c) = e(g_1)$. The case of $\alpha = c$ is justified as this one. \\

In particular, $\mathcal{K} \lhd K$ and since $\mathcal{K}$ contains $K'$, the quotient $K/\mathcal{K}$ is abelian and generated by the classes of $x,y,z,t$ and $w$. Notice that $yt^{-1} = (xy^{-1},1)$ and $zw^{-1} = (1, xy^{-1})$. Hence, $[y] = [t]$ and $[z] = [w]$ modulo $\mathcal{K}$, which implies that the classes of $x,y,z$ in fact generate $K/\mathcal{K}$. Furthermore, $y^2 = (x^2,1)$ and $z^2 = (1,x^2)$, so the orders of $[y]$ and $[z]$ are $2$, proving that $K/\mathcal{K} \leq C_4 \times C_2 \times C_2$. It remains to prove that there are no more relations. 

Let $g = x^\alpha y^\beta z^\gamma$ with $\alpha \in \set{0,\dots,3}$ and $\beta, \gamma \in \set{0,1}$, and suppose that $g$ represents the trivial class modulo $\mathcal{K}$. First, we have that $x^\alpha = ((ca)^\alpha, (ac)^\alpha)$. The classes of the elements $ac$ and $ca$ have order $4$ in $G/K$ and therefore $x^\alpha \in K_1$ if and only if $\alpha = 0$. Since $y$ and $z \in K_1$, and we are assuming that $g \in \mathcal{K} \leq K_1$, necessarily $\alpha$ must be zero. Thus, $g = y^\beta z^\gamma = (x^\beta, x^\gamma)$ and since $g \in \mathcal{K}$, necessarily $\beta = \gamma = 0$. \\

Finally, by \cref{theorem: K/K' is C45}, $[K:K'] = 4^5$ and $[K:\mathcal{K}] = 4^2$, so
\begin{equation*}
[\mathcal{K}: K'] = 4^3 = 64. \qedhere
\end{equation*}
\end{proof}

\section{Branch kernel and rigid kernel}
\label{section: branch kernel and rigid kernel}

We finish this article by computing the branch kernel and the rigid kernel of $G$.

\begin{Theorem}
The rigid kernel of $G$ is trivial and the branch kernel is isomorphic to $C_4[\partial T]$.
\end{Theorem}

\begin{proof}
By \cref{lemma: equivalences inclusion profinite kernels}, the rigid kernel is trivial if and only if for all $n \in \NN$ there exists $m \in \NN$ such that $\RiSt_G(n) \geq \St_G(m)$. To verify this, observe that $K_n$ is contained in $\RiSt_G(n)$. Since $K \geq \St_G(3)$, then $$\RiSt_G(n) \geq K_n \geq (\St_G(3))_n \simeq \St_G(n+3),$$ where in the last step \cref{lemma: St(n+1) = prod St(n)} was used.

We now prove the result for the branch kernel. By \cite[Remark 2.8]{Bartholdi2012}, the branch kernel is isomorphic to $A[[\partial T]]$ and by \cite[Remark 2.10]{Bartholdi2012} the group $A$ is calculated in the following way: consider the map $\sigma: K \rightarrow K$ defined as $\sigma(g) = (1,g)$. The map is well defined since $K$ is regular branch. Moreover, since $K' \geq (K')_1$, the map descends to a well-defined map $\sigma: K/K' \rightarrow K/K'$. Then, the subgroup $A$ is defined as $$A = \bigcap_{n \geq 0} \sigma^n(K/K').$$ Notice that the sequence of subgroups considered in the intersection is nested, so we need to determine when it stabilizes.

Since $K/K'$ is abelian, we will use additive notation. We start calculating the images of $xK',\dots,wK'$. To simplify the notation even more, we will simply write $x,\dots,w$ for the images of them in $K/K'$. For the calculation, we will use \cref{lemma: relations conjugations K IMG} many times: 

\begin{itemize}
\item $\sigma(x) = (1,x) = z$ 
\item $\sigma(y) = (1,y) = w$
\item $\sigma(z) = (1,z) = (1,aya) = (aba)w(aba) = (ab)t(ba) = a(-t)a = -w$
\item $\sigma(t) = (1,t) = (1,cxcy^{-1}x^{-1}) = bzb -w - z = -w -2z$
\item $\sigma(w) = (1,w) = (1,ata) = aba(-2z-w)aba = ab(-2y -t)ba = a(2y - t)a = 2z-w$
\end{itemize}

Hence, $\sigma(K/K') = \< z,w,-w+2z,-w-2z,w > = \<z,w>$. Then, we just iterate the previous calculations:
$$\sigma^2(K/K') = \sigma(\<z,w>) = \<2z-w, -w> = \<2z+w, w>.$$

Next, observe that $\sigma(2z+w) = 2z+w$ and $\sigma(w) = 2z-w$, so $$\sigma^3(K/K') = \sigma(\<2z+w,w>) = \<2z+w, 2z-w> = \<2z+w, 2w>,$$ but $2w = 2(2z+w)$ because $z$ has order $4$. Thus, $$\sigma^3(K/K') = \<2z+w>$$ and since it is fixed by the applications of $\sigma$, we have $\sigma^n(K/K') = \<2z+w>$ for all $n \geq 3$. Therefore, $A = \<2z+w> \simeq C_4$.
\end{proof}

\bibliographystyle{unsrt}

\end{document}